\documentclass{article}
\usepackage{amsmath,amssymb,amsfonts,amsthm}
\usepackage{graphicx}
\usepackage{layout}
\usepackage[latin1]{inputenc}
\usepackage[french,english]{babel}
\usepackage[french]{babel}
\usepackage{cancel}

\newtheorem{theorem}{Theorem}

\newtheorem{corollary}{Corollary}
\newtheorem{proposition}{Proposition}
\newtheorem{lemma}{Lemma}

\newtheorem{remark}{Remark}    

\usepackage{fancyhdr}

\usepackage{setspace}

\renewcommand{\headrulewidth}{0,5pt}
\newcommand{\vertiii}[1]{{\left\vert\kern-0.25ex\left\vert\kern-0.25ex\left\vert #1 
    \right\vert\kern-0.25ex\right\vert\kern-0.25ex\right\vert}} 

\title{Perturbation of Schauder frames and  besselian Schauder frames
 of  Banach spaces}
\author{{\footnotesize Samir Kabbaj, Rafik Karkri and Hicham Zoubeir }  }

\date{\today}

\begin{document}
\maketitle 
\begin{abstract}
 We consider the stability of Schauder frames and
besselian Schauder frames under perturbations. Our results are 
inspirit close to the results of Heil \cite{C.Heil.1990}. 
  
\end{abstract}
{\footnotesize $\textbf{MSC2020-Mathematics Subject Classification System}$. 40D05, 46B03,46B15,
46B20, 46B25, 46B45, 47A30, 47A55, 47A63, 47B01, 47B38.}\\
{\footnotesize $\textbf{Keywords and phrases:}$ Separable Banach space, Weakly 
sequentially complete Banach space, Besselian sequence,
Frame, Schauder frame, Besselian schauder frame, Perturbation, 
Perturbation of Schauder frames.}

\renewcommand{\headrulewidth}{0pt}  
\fancyhead[C]{{\scriptsize Samir Kabbaj, Rafik Karkri  and Hicham Zoubeir }} 

\section{Introduction}
In $1946,$ Gabor \cite{gab} performed a new method for the signal
reconstruction from elementary signals. In $1952,$
 Duffin and Schaeffer \cite{duf} developped, in the field of nonharmonic series, a
similar tool and introduced frame theory for Hilbert spaces.
A frames for a Banach space were 
 introduced in $(1991)$ by Gr\"{o}chenig \cite{gro}. 
 In light of the works of  Cassaza, Han and Larson  \cite
{lar1},    Han and Larson \cite{lar2}.
 Cassaza \cite{cas.2008}   introduced the notion of
Schauder frames of a given  Banach spaces. In 2021 Karkri and Zoubeir
  \cite{karkri.zoubeir.2021} introduced the notion of
besselian Schauder frames of Banach spaces and obtained some results  on
universal  spaces and complemented subspaces for Banach spaces with Schauder frames 
or with besselian Schauder frames.
 For more generalizations  of  frame notion one can refer to 
 \cite{cas2, Naroei.Nazari, Kim.Lim, Hua.Yongdong}.

The first perturbation result of Schauder basis in Banach spaces
 is due to Paley and Wiener \cite{Paley.Wiener.1934}. 
 Frames perturbation of atoms in Banach spaces were introduced by Heil in his Ph.D.thesis 
\cite{C.Heil.1990}. The result of Heil  and the book 
\cite{R.Young.1980} motivate O.Christensen
\cite{Christensen.1995}  to consider a similar  problem
for frames in Hilbert spaces. In 1997 Casazza and Christensen 
\cite{Casazza.Christensen.1997} proved new theorems of perturbation
of frames in both Hilbert spaces and Banach spaces.
In this paper we obtain new theorems of perturbation of Schauder
frames and besselian Schauder frames of Banach spaces.
\section{Main definitions and notations}
Let $X$  be a separable Banach space on $\mathbb{K}\in \left \{ \mathbb{R},
\mathbb{C}\right \} $ and  $X^{\ast }$ it topological dual. 
\begin{enumerate}
\item
The Banach space $X$ is said to be weakly sequentially complete  if for each
sequence $(x_{n})_{n\in \mathbb{N}^{* }}$ of $X$ such that 
$\underset{n\rightarrow +\infty }{\lim }x^{* }(x_{n})$ exists for every $x^{*
}\in X^{* },$ there exists $x\in X$ such that $\underset{n\rightarrow
+\infty }{\lim }x^{* }(x_{n})=x^{* }(x)$ for every $x^{* }\in
X^{* }$ \cite[page 218, definition 2.5.23]{Megginson}
 \cite[pages 37-38]{kal}.
\item We denote by $\mathbb{B}_{X\text{ \ }}$ the closed unit ball of $X:$%
\begin{equation*}
\mathbb{B}_{X}:=\left \{ x\in X:\left \Vert x\right \Vert _{X}\leq 1\right \}
\end{equation*}
\item We denote by $L\left( X\right) $ the set of all\ bounded linear
operators $f:X\rightarrow X$ . It is well-known that $L\left( X\right) $ is
a Banach space for the norm $\left \Vert \cdot \right \Vert _{L\left(
X\right) }$ defined by the formula:
\begin{equation*}
\left \Vert f\right \Vert _{L\left( X\right) }:=\underset{x\in \mathbb{B}_{X}%
}{\sup }\left \Vert f\left( x\right) \right \Vert _{X}
\end{equation*}
\item
For a sequence $\left(  x_{n}\right)_{n\in\mathbb{N}^{*}}$ of elements  of 
$X$ we denote by $span \left(  x_{n}\right)_{n\in\mathbb{N}^{*}}$ the algebraic 
linear span of $\left(  x_{n}\right)_{n\in\mathbb{N}^{*}}$ i.e. the set of all 
finite linear combinations of $\left(  x_{n}\right)_{n\in\mathbb{N}^{*}}$.
The closure of $span \left(  x_{n}\right)_{n\in\mathbb{N}^{*}}$ is denoted by
$\underset{n\in \mathbb{N}^{*}}{\overline{span}}\left (x_{n} \right )$.
\item A sequence 
$\mathcal{F}=\left( \left(x_{n},y_{n}^{* }\right) \right)_{n\in\mathbb{N}^{*}}\subset X\times X^{*} $ 
is called a paire of $X$ \cite{karkri.zoubeir.2021}.
\item The paire   $\mathcal{F}$ of $X$ is called a  Schauder frame 
(resp. unconditional  Schauder frame ) of $X$ if for all $x\in X$, 
the series $\sum y_{n}^{\ast }\left( x\right) x_{n}$
 is convergent (resp. unconditionally convergent) in $X$ 
 to $x$ \cite{cas.2008}, \cite{karkri.zoubeir.2021}.
\item
If $\mathcal{F}$ is a Schauder frame of $X$, we denote by $K_{\mathcal{F}}$ the finite quantity
 $$ K_{\mathcal{F}}:=\underset{x\in \mathbb{B}_{E},\text{ }n\in\mathbb{N}^{\ast }}{\sup }\left( \left \Vert \underset{k=1}{\overset{n}{\sum }}
b_{k}^{\ast }\left( x\right) a_{k}\right \Vert \right) $$
\item   $\left( \left( x_{n},y_{n}^{\ast }\right) \right) _{n\in 
\mathbb{N}^{*}}$   is said to be a besselian paire of $X$ if there exists a
constant $A>0$ such that 
\begin{equation*}
\underset{k=1}{\overset{+\infty}{\sum }}\left \vert y_{n}^{\ast }\left( x\right)
\right \vert \left \vert y^{\ast }\left( x_{n}\right) \right \vert \leq
A\left \Vert x\right \Vert _{X}\left \Vert y^{\ast }\right \Vert _{X^{\ast }}
\end{equation*}%
for each $\left( x,y^{\ast }\right) \in X\times X^{\ast }$ 
\cite{karkri.zoubeir.2021}.
\item $\left( \left( x_{n},y_{n}^{\ast }\right) \right) _{n\in 
\mathbb{N}^{*}}$   is said to be a besselian Schauder frame of $X$ if it is
both a Schauder frame and a besselian paire \cite{karkri.zoubeir.2021}.
\item
Let $\left( x_{n}\right) _{n\in\mathbb{N}^{* }}$ and 
$\left( w_{n}\right) _{n\in\mathbb{N}^{* }}$ be sequences of 
elements of $X$. We use the notation
$ \left(w_{n} \right)_{n\in\mathbb{N}^{*}}\equiv_{X}
 \left(x_{n} \right)_{n\in\mathbb{N}^{*}}  $
 if there is an isomorphism $T$ from $X$ onto $X$ such that 
$w_{n}=T(x_{n})$ for each $n\in \mathbb{N}^{*}$. 
\item
Let $\left( y^{*}_{n}\right) _{n\in\mathbb{N}^{* }}$ and 
$\left( z^{*}_{n}\right) _{n\in\mathbb{N}^{* }}$ be sequences of 
elements of $X^{*}$. We use the notation
$ \left(z^{*}_{n} \right)_{n\in\mathbb{N}^{*}}\equiv_{X^{*}}
 \left(y^{*}_{n} \right)_{n\in\mathbb{N}^{*}}  $
 if there is an isomorphism $T$ from $X$ onto $X$ such that 
$z^{*}_{n}=y^{*}_{n}\circ T$ for each $n\in \mathbb{N}^{*}$. 
\end{enumerate}
 For all the material on Banach spaces, Hilbertian frames and
 more generalizations  of  frame notion   to
 Banach spaces , one can refer to
\cite{Megginson}, \cite{lin01}, \cite{lin02}, \cite{woj}, \cite{cas2}, 
\cite{fei}, \cite{chr3}, \cite{cas1}, \cite{dau},
 \cite{lar2}, \cite{R.Young.1980}, \cite{ald}, 
  \cite{chr1} and \cite{Zizler.2011}.
  In the sequel $E$ will represent a given separable Banach space and 
  $\mathcal{F}=\left( \left(a_{n},b_{n}^{* }\right) \right) _{n\in\mathbb{N}^{*}} $
a paire of $E$.
     
   \section{Fundamental results}
\begin{theorem}
\label{thm.pert.Schauder.f}
Let $\mathcal{F}=\left( \left(a_{n},b_{n}^{* }\right) \right) _{n\in\mathbb{N}^{*}} $
 be a Schauder frame of $E$ and 
 $\left( \left(x_{n},y_{n}^{* }\right) \right) _{n\in\mathbb{N}^{*}}$ a paire of
 $E$. We assume that:
\begin{align}
\label{perturbation.condition}
\underset{n= 1}{\overset{\infty}{\sum }} \left \Vert y^{*}_{n}-b^{*}_{n}\right \Vert_{E^{*}}\left \Vert x_{n}\right \Vert_{E}
+\underset{a_{n}\neq 0}{\overset{}{\sum }}2K_{\mathcal{F}}\dfrac{ \left\Vert x_{n}-a_{n}\right \Vert_{E}}{\left \Vert a_{n}\right \Vert_{E}} 
+\underset{a_{n}= 0}{\overset{}{\sum }} \left \Vert b^{*}_{n}\right \Vert_{E^{*}}\left \Vert x_{n}\right \Vert_{E}
<1
\end{align}
Then there exists a sequence $ \left(z_{n}^{*}\right)_{n\in\mathbb{N}^{*}} $  
 of elements of $E^{*}$ and a sequence $ \left(w_{n}\right)_{n\in\mathbb{N}^{*}} $  
  of elements of $E$ such that the following conditions holds:
\begin{enumerate}
\item $ \left(w_{n} \right)_{n\in\mathbb{N}^{*}}\equiv_{E}
 \left(x_{n} \right)_{n\in\mathbb{N}^{*}}  $. 
 \item $ \left(z_{n}^{*} \right)_{n\in\mathbb{N}^{*}}\equiv_{E^{*}}
 \left(y_{n}^{*} \right)_{n\in\mathbb{N}^{*}}$.
\item    $\left( \left(x_{n},z_{n}^{* }\right) \right) _{n\in\mathbb{N}^{*}} $
and $\left( \left(w_{n},y_{n}^{* }\right) \right) _{n\in\mathbb{N}^{*}} $
 are a Schauder frames of $E$.
\end{enumerate}
\end{theorem}
$ \textbf{Proof} $: 
Let $x\in E$ and  $n\in\mathbb{N}^{*}$. We have:
\[ y_{n}^{*}(x)x_{n}-b_{n}^{*}(x)a_{n}=(y_{n}^{*}-b_{n}^{*})(x)x_{n}+b_{n}^{*}(x)(x_{n}-a_{n})\]
It follows that:
\[\left\Vert y_{n}^{*}(x)x_{n}-b_{n}^{*}(x)a_{n}\right \Vert_{E}
\leq  \left\Vert y_{n}^{*}-b_{n}^{*}\right \Vert_{E^{*}}\left\Vert x_{n}\right \Vert_{E} \left\Vert x\right \Vert_{E}
+\left\vert b_{n}^{*}(x)\right \vert \left\Vert x_{n}-a_{n}\right \Vert_{E}\]
If $a_{n}\neq 0$, then:
\begin{align*}
 \left\vert b_{n}^{*}(x)\right \vert\left\Vert   x_{n}-a_{n}\right \Vert_{E}&=\left\Vert b_{n}^{*}(x)a_{n}\right \Vert_{E} \dfrac{\left\Vert   x_{n}-a_{n}\right \Vert_{E}}{\left\Vert  a_{n}\right \Vert_{E}}\\
&\leqslant 2K_{\mathcal{F}}\left\Vert x\right \Vert_{E} \dfrac{\left\Vert   x_{n}-a_{n}\right \Vert_{E}}{\left\Vert   a_{n}\right \Vert_{E}}
\end{align*}
Hence: 
\[\left\Vert y_{n}^{*}(x)x_{n}-b_{n}^{*}(x)a_{n}\right \Vert_{E}
\leq \left(  \left\Vert y_{n}^{*}-b_{n}^{*}\right \Vert_{E^{*}}\left\Vert x_{n}\right \Vert_{E} 
+2K_{\mathcal{F}} \dfrac{\left\Vert   x_{n}-a_{n}\right \Vert_{E}}{\left\Vert   a_{n}\right \Vert_{E}}\right) \left\Vert x\right\Vert_{E} \]
If  $a_{n}=0$. Then,
\[\left\vert b_{n}^{*}(x)\right\vert \left\Vert x_{n}-a_{n}\right \Vert_{E}
 \leq \left\Vert b_{n}^{*}\right \Vert_{E^{*}} \left\Vert x\right \Vert_{E} \left\Vert x_{n}\right \Vert_{E}
\]
Hence:
\[\left\Vert y_{n}^{*}(x)x_{n}-b_{n}^{*}(x)a_{n}\right \Vert_{E}
\leq \left(  \left\Vert y_{n}^{*}-b_{n}^{*}\right \Vert_{E^{*}}\left\Vert x_{n}\right \Vert_{E} 
+\left\Vert b_{n}^{*}\right \Vert_{E^{*}}\left\Vert x_{n}\right \Vert_{E} \right) \left\Vert x\right\Vert_{E} \]
Consequently, we have for each $x\in E$ and $n\in \mathbb{N}^{*}$
$$\underset{n=1}{\overset{+\infty}{\sum }}\left\Vert y_{n}^{*}(x)x_{n}-b_{n}^{*}(x)a_{n}\right \Vert_{E} \leq Q \left\Vert x\right \Vert_{E}$$
where
$$Q:=\underset{n=1}{\overset{+\infty}{\sum }} \left\Vert y_{n}^{*}- b_{n}^{*}\right \Vert_{E^{*}}  \left\Vert x_{n}\right \Vert_{E}
+ \underset{a_{n}\neq 0}{\overset{}{\sum }} 2K_{\mathcal{F}} \dfrac{\left\Vert   x_{n}-a_{n}\right \Vert_{E}}{\left\Vert   a_{n}\right \Vert_{E}}
+\underset{a_{n}= 0}{\overset{}{\sum }} \left\Vert b_{n}^{*}\right \Vert_{E^{*}}  \left\Vert x_{n}\right \Vert_{E}$$
We have according to the inequality (\ref{perturbation.condition}) then $Q\in [0,1[$. It 
follows that the series $\sum_{n}y_{n}^{*}(x)x_{n}$ is convergent in $E$ for each 
$x\in E$. We set
\[T(x)=\underset{n=1}{\overset{+\infty}{\sum }}y_{n}^{*}(x)x_{n}, \;\; x\in E\]
We have then 
\begin{align*}
 				 \left\Vert T(x)-x\right \Vert_{E}& =\left\Vert \underset{n=1}{\overset{+\infty}{\sum }}y_{n}^{*}(x)x_{n}
 				 -\underset{n=1}{\overset{+\infty}{\sum }}b_{n}^{*}(x)a_{n}\right \Vert_{E}\\	  
 				 	 &\leq \underset{n=1}{\overset{+\infty}{\sum }}\left\Vert y_{n}^{*}(x)x_{n}-b_{n}^{*}(x)a_{n}\right \Vert_{E}\\
 				 	 & \leq Q \left\Vert x\right \Vert_{E}
\end{align*}
It follows that $T: E\longrightarrow E$ is an isomorphism of Banach spaces. Consequently,
\[x=\underset{n=1}{\overset{+\infty}{\sum }}y_{n}^{*}(x)T^{-1}(x_{n}), \;\; x\in E\]
and
\[x=\underset{n=1}{\overset{+\infty}{\sum }}(y_{n}^{*}\circ T^{-1})(x)x_{n}, \;\; x\in E\]
Let us set 
\[\left\{
\begin{array}{l}
 w_{n}=T^{-1}(x_{n}), \;\; x\in E \\
 z_{n}^{*}=y_{n}^{*}\circ T^{-1}, \;\; x\in E
\end{array}
\right.\]
Hence
\[\left\{
\begin{array}{l}
\left(w_{n} \right)_{n\in\mathbb{N}^{*}}\equiv_{E}
 \left(x_{n} \right)_{n\in\mathbb{N}^{*}} \\
 \left(z_{n}^{*} \right)_{n\in\mathbb{N}^{*}}\equiv_{E^{*}}
 \left(y_{n}^{*} \right)_{n\in\mathbb{N}^{*}}
\end{array}
\right.\]
Furthermore, 
$\left( \left(x_{n},z_{n}^{* }\right) \right) _{n\in\mathbb{N}^{*}} $
and $\left( \left(w_{n},y_{n}^{* }\right) \right) _{n\in\mathbb{N}^{*}} $
 are a Schauder frames of $E$.\qed
\begin{lemma}\cite{karkri.zoubeir.2021}
\label{unconditionally}
We assume that  $E$\textit{\ is a weakly
sequentially complete Banach space and that} $\mathcal{F}:=\left( \left(
a_{n},b_{n}^{* }\right) \right) _{n\in 
\mathbb{N}
^{* }}$ \textit{is a besselian paire of }$E.$ \textit{Then for each }$x\in
E$\textit{, the series }$\sum b_{n}^{* }\left( x\right) a_{n}$\textit{\
is unconditionally convergent in }$E.$
\end{lemma}
\noindent {\bf Proof.}
For each $x\in E$, $y^{* }\in E^{* }$ we have: 
\begin{align*}
\overset{+\infty }{\underset{n=1}{\sum }}\left \vert y^{* }\left( b_{n}^{*
}\left( x\right) a_{n}\right) \right \vert & =\overset{+\infty }{\underset{%
k=1}{\sum }}\left \vert b_{n}^{* }\left( x\right) y^{* }\left( a_{n}\right)
\right \vert \\
& \leq \mathcal{L}_{\mathcal{F}}\left \Vert x\right \Vert _{E}\left \Vert
y^{* }\right \Vert _{E^{* }} \\
& <+\infty
\end{align*}
Hence the series $\sum b_{n}^{* }\left( x\right) a_{n}$ is weakly
unconditionally convergent. Then, since $E$ is weakly sequentially complete,
the well-known Orlicz's theorem $(1929)$ \cite[Proposition.4 ,page 59 and page 66]{woj}, entails that the series $\sum b_{n}^{* }\left( x\right) a_{n}$
is unconditionally   convergent.
$\Box$
\begin{theorem}
\label{thm.pert.b.Schauder.f}
Let $E$ be a weakly sequentially complete Banach space and 
 $\mathcal{F}=\left( \left(a_{n},b_{n}^{* }\right) \right) _{n\in\mathbb{N}^{*}} $
 be a besselian Schauder frame of $E$. Let $\left( \left(x_{n},y_{n}^{* }\right) \right) _{n\in\mathbb{N}^{*}} $
 be a paire of $E$ for which there exists a constant $\alpha\in [0,1[$
 such that the following condition holds for each $x\in E$ and 
 $y^{*}\in E^{*}$
\begin{align*}
\underset{n= 1}{\overset{\infty}{\sum }} \left \vert y^{*}_{n}(x)y^{*}(x_{n}-a_{n})+(y^{*}_{n}-b^{*}_{n})(x)y^{*}(a_{n})\right \vert
<\alpha \left \Vert x\right \Vert_{E}\left \Vert y^{*}\right \Vert_{E^{*}}
\end{align*}
Then there exists a sequence $ \left(z_{n}^{*}\right)_{n\in\mathbb{N}^{*}} $  
 of elements of $E^{*}$ and a sequence $ \left(w_{n}\right)_{n\in\mathbb{N}^{*}} $  
  of elements of $E$ such that the following conditions holds:
\begin{enumerate}
\item $ \left(w_{n} \right)_{n\in\mathbb{N}^{*}}\equiv_{E}
 \left(x_{n} \right)_{n\in\mathbb{N}^{*}}  $. 
 \item $ \left(z_{n}^{*} \right)_{n\in\mathbb{N}^{*}}\equiv_{E^{*}}
 \left(y_{n}^{*} \right)_{n\in\mathbb{N}^{*}}$.
\item    $\left( \left(x_{n},z_{n}^{* }\right) \right) _{n\in\mathbb{N}^{*}} $
and $\left( \left(w_{n},y_{n}^{* }\right) \right) _{n\in\mathbb{N}^{*}} $
 are a besselian Schauder frames of $E$.
\end{enumerate}
 \end{theorem}
$\textbf{Proof}$: Let $x\in E$ and $n\in\mathbb{N}^{*}$. We have:
\[y^{*}_{n}(x)x_{n}-b^{*}_{n}(x)a_{n}
=y^{*}_{n}(x)(x_{n}-a_{n})
+ (y^{*}_{n}-b^{*}_{n})(x)a_{n} \]
It follows that we have for all $x\in E$ and $y^{*}\in E^{*}$:  
\begin{align*}
 				 	 \underset{n= 1}{\overset{\infty}{\sum }} \left \vert y^{*}_{n}(x)y^{*}(x_{n}) \right \vert
 				 	 & \leq  \underset{n= 1}{\overset{\infty}{\sum }} \left \vert b^{*}_{n}(x)y^{*}(a_{n}) \right \vert
 				 	 + \underset{n= 1}{\overset{\infty}{\sum }} \left \vert y^{*}_{n}(x)y^{*}(x_{n}-a_{n})+(y^{*}_{n}-b^{*}_{n})(x)y^{*}(a_{n}) \right \vert\\
 				   	 &\leq (\alpha+ \mathcal{L}_{\mathcal{F}})\left \Vert x \right \Vert_{E}\left \Vert y^{*} \right \Vert_{E^{*}}
\end{align*}
Since $E$ is a weakly sequentially complet Banach space,
 it follows from lemma \ref{unconditionally} that the series
$\sum _{n}y^{*}_{n}(x)x_{n}$ is convergent.
 Let us then set:
 \[T(x)=\underset{n= 1}{\overset{\infty}{\sum }}y^{*}_{n}(x)x_{n} \]
Then we have:
\begin{align*}
 				 	  \left \Vert T(x)-x \right \Vert_{E}
 				 	 & =\left \Vert\underset{n= 1}{\overset{\infty}{\sum }}  y^{*}_{n}(x)(x_{n}-a_{n})+(y^{*}_{n}-b^{*}_{n})(x)a_{n} \right \Vert_{E}\\
 	&=\underset{y^{*}\in \mathbb{B}_{E^{*}}}{\sup}\left \vert\underset{n= 1}{\overset{\infty}{\sum }}  y^{*}_{n}(x)y^{*}(x_{n}-a_{n})+(y^{*}_{n}-b^{*}_{n})(x)y^{*}(a_{n}) \right \vert\\			   	 
 				   	 &\leq \alpha\left \Vert x \right \Vert_{E}
\end{align*}
It follows that $T$ is a bounded linear operator on $E$  and that:
\[\left \Vert T-Id_{E}\right \Vert_{L(E)}\leq \alpha <1\] 
Hence $T$ is an isomorphism of Banach spaces. Consequently, the following
identities holds for each $x\in E$:
\[\left\{
\begin{array}{l}
 x=\underset{n=1}{\overset{+\infty}{\sum }}(y_{n}^{*}\circ T^{-1})(x)x_{n}  \\
x=\underset{n=1}{\overset{+\infty}{\sum }}y_{n}^{*} (x)T^{-1}(x_{n}) 
\end{array}
\right.\]
Let us then set 
\[\left\{
\begin{array}{l}
 w_{n}=T^{-1}(x_{n}),\;\;n\in \mathbb{N}^{*} \\
z^{*}_{n}=y_{n}^{*}\circ T^{-1},\;\;n\in \mathbb{N}^{*}
\end{array}
\right.\]
Hence
\[\left\{
\begin{array}{l}
 \left(w_{n} \right)_{n\in\mathbb{N}^{*}}\equiv_{E}
 \left(x_{n} \right)_{n\in\mathbb{N}^{*}}  \\
 \left(z_{n}^{*} \right)_{n\in\mathbb{N}^{*}}\equiv_{E^{*}}
 \left(y_{n}^{*} \right)_{n\in\mathbb{N}^{*}}
\end{array}
\right.\]
Furthermore, we have for each $x\in E$ and $y^{*}\in E^{*}$
\begin{align*}
 				 	 \underset{n= 1}{\overset{\infty}{\sum }} \left \vert y^{*}_{n}(x)y^{*}(w_{n}) \right \vert
 				 	 & =  \underset{n= 1}{\overset{\infty}{\sum }} \left \vert y^{*}_{n}(x)y^{*}(T^{-1}(x_{n})) \right \vert \\
 				 	  & =  \underset{n= 1}{\overset{\infty}{\sum }} \left \vert y^{*}_{n}(x)(y^{*}\circ T^{-1})(x_{n}) \right \vert \\
 				 	 &\leq (\alpha+ \mathcal{L}_{\mathcal{F}})\left \Vert x \right \Vert_{E}\left \Vert y^{*}\circ T^{-1} \right \Vert_{E^{*}} \\
	 &\leq (\alpha+ \mathcal{L}_{\mathcal{F}})\left \Vert T^{-1}\right \Vert_{L(E)}\left \Vert x \right \Vert_{E}\left \Vert y^{*} \right \Vert_{E^{*}}
\end{align*}
and
\begin{align*}
 				 	 \underset{n= 1}{\overset{\infty}{\sum }} \left \vert z^{*}_{n}(x)y^{*}(x_{n}) \right \vert
 				 	 & =  \underset{n= 1}{\overset{\infty}{\sum }} \left \vert y^{*}_{n}(T^{-1}(x))y^{*}(x_{n}) \right \vert \\
 				 	 &\leq (\alpha+ \mathcal{L}_{\mathcal{F}})\left \Vert T^{-1}(x) \right \Vert_{E}\left \Vert y^{*}  \right \Vert_{E^{*}} \\
	 &\leq (\alpha+ \mathcal{L}_{\mathcal{F}})\left \Vert T^{-1}\right \Vert_{L(E)}\left \Vert x \right \Vert_{E}\left \Vert y^{*} \right \Vert_{E^{*}}
\end{align*}
It follows that  $\left( \left(x_{n},z_{n}^{* }\right) \right) _{n\in\mathbb{N}^{*}} $
and $\left( \left(w_{n},y_{n}^{* }\right) \right) _{n\in\mathbb{N}^{*}} $
 are a besselian Schauder frames of $E$.\qed
\begin{corollary}
\label{pert.bsf}
Let $E$ be a weakly sequentially complete Banach space and 
$\left( \left(a_{n},b_{n}^{* }\right) \right) _{n\in\mathbb{N}^{*}} $ be a besselian Schauder frame of $E$. 
Let $\left( \left(x_{n},y_{n}^{* }\right) \right) _{n\in\mathbb{N}^{*}} $
 be a paire of $E $ such that 
 \begin{align*}
\underset{n= 1}{\overset{\infty}{\sum }} \left \Vert y^{*}_{n}\right \Vert_{E^{*}} \left \Vert x_{n}-a_{n} \right \Vert_{E}
+\left \Vert y^{*}_{n}-b^{*}_{n} \right \Vert_{E^{*}} \left \Vert a_{n}\right \Vert_{E}
<1
\end{align*} 
Then there exists a sequence $ \left(z_{n}^{*}\right)_{n\in\mathbb{N}^{*}} $  
 of elements of $E^{*}$ and a sequence $ \left(w_{n}\right)_{n\in\mathbb{N}^{*}} $  
  of elements of $E$ such that the following conditions holds:
\begin{enumerate}
\item $ \left(w_{n} \right)_{n\in\mathbb{N}^{*}}\equiv_{E}
 \left(x_{n} \right)_{n\in\mathbb{N}^{*}}  $. 
 \item $ \left(z_{n}^{*} \right)_{n\in\mathbb{N}^{*}}\equiv_{E^{*}}
 \left(y_{n}^{*} \right)_{n\in\mathbb{N}^{*}}$.
\item    $\left( \left(x_{n},z_{n}^{* }\right) \right) _{n\in\mathbb{N}^{*}} $
and $\left( \left(w_{n},y_{n}^{* }\right) \right) _{n\in\mathbb{N}^{*}} $
 are a besselian Schauder frames of $E$.
\end{enumerate}
\end{corollary}
$\textbf{Proof}$: Let us set: $\alpha:=\underset{n= 1}{\overset{\infty}{\sum }} \left \Vert y^{*}_{n}\right \Vert_{E^{*}} \left \Vert x_{n}-a_{n} \right \Vert_{E}
+\left \Vert y^{*}_{n}-b^{*}_{n} \right \Vert_{E^{*}} \left \Vert a_{n}\right \Vert_{E}
$. We have then $\alpha\in [0,1[$ and for each $x\in E$ and $y^{*}\in E^{*}$:
\begin{align*}
\underset{n= 1}{\overset{\infty}{\sum }} \left \vert y^{*}_{n}(x)y^{*}(x_{n}-a_{n})+(y^{*}_{n}-b^{*}_{n})(x)y^{*}(a_{n})\right \vert
<\alpha \left \Vert x\right \Vert_{E}\left \Vert y^{*}\right \Vert_{E^{*}}
\end{align*}
Consequently, by virtue of the theorem \ref{thm.pert.b.Schauder.f}, there exists
  a sequence $ \left(z_{n}^{*}\right)_{n\in\mathbb{N}^{*}} $
 of elements of  $E^{*}$ and  
a sequence $ \left(w_{n} \right)_{n\in\mathbb{N}^{*}} $  of elements of 
 $E$ such that $\left( \left(x_{n},z_{n}^{* }\right) \right) _{n\in\mathbb{N}^{*}} $
and $\left( \left(w_{n},y_{n}^{* }\right) \right) _{n\in\mathbb{N}^{*}} $
 are a besselian Schauder frames of $E$.\qed
\begin{corollary}
\label{pertubation.1}
Let $E$ be a weakly sequentially complete Banach space and 
$\left( \left(a_{n},b_{n}^{* }\right) \right) _{n\in\mathbb{N}^{*}} $ 
be a besselian Schauder frame of $E$. 
Let $\left( x_{n} \right) _{n\in\mathbb{N}^{*}} $
 be a sequence of elements of $E $ for which there exists a constant
  $\alpha\in [0,1[$ such that the following condition holds for each 
  $x\in E$ and $y^{*} \in E^{*}$:
 \begin{align*}
\underset{n= 1}{\overset{\infty}{\sum }} \left \vert b^{*}_{n}(x)y^{*}(x_{n}-a_{n})\right \vert 
\leq \alpha \left\Vert x\right \Vert_{E}\left\Vert y^{*}\right \Vert_{E^{*}}
\end{align*} 
Then there exists a sequence $ \left(z_{n}^{*}\right)_{n\in\mathbb{N}^{*}} $  
 of elements of $E^{*}$ and a sequence $ \left(w_{n}\right)_{n\in\mathbb{N}^{*}} $  
  of elements of $E$ such that the following conditions holds:
\begin{enumerate}
\item $ \left(w_{n} \right)_{n\in\mathbb{N}^{*}}\equiv_{E}
 \left(x_{n} \right)_{n\in\mathbb{N}^{*}}  $. 
 \item $ \left(z_{n}^{*} \right)_{n\in\mathbb{N}^{*}}\equiv_{E^{*}}
 \left(b_{n}^{*} \right)_{n\in\mathbb{N}^{*}}$.
\item    $\left( \left(x_{n},z_{n}^{* }\right) \right) _{n\in\mathbb{N}^{*}} $
and $\left( \left(w_{n},y_{n}^{* }\right) \right) _{n\in\mathbb{N}^{*}} $
 are a besselian Schauder frames of $E$.
\end{enumerate}
\end{corollary}
$\textbf{Proof}$: Let us consider the paire 
$\left( \left(x_{n},b_{n}^{* }\right) \right) _{n\in\mathbb{N}^{*}} $  of $E$.
 We have then for each $x\in E$ and $y^{*}\in E^{*}$:
 \begin{align*}
\underset{n= 1}{\overset{\infty}{\sum }} \left \vert b^{*}_{n}(x)y^{*}(x_{n}-a_{n})+(b^{*}_{n}-b^{*}_{n})(x)y^{*}(a_{n})\right \vert 
\leq \alpha \left\Vert x\right \Vert_{E}\left\Vert y^{*}\right \Vert_{E^{*}}
\end{align*} 
It follows then, by virtue of the theorem \ref{thm.pert.b.Schauder.f} , that there exists a sequence
$ \left(z_{n}^{*}\right)_{n\in\mathbb{N}^{*}} $ of elements of $E^{*}$ and a sequence
$ \left(w_{n}\right)_{n\in\mathbb{N}^{*}} $ of elements of $E$ such that 
the conditions (1), (2) and (3) in the corollary are true. \qed
\begin{corollary}
\label{pertubation.2}
Let $E$ be a weakly sequentially complete Banach space and 
$\left( \left(a_{n},b_{n}^{* }\right) \right) _{n\in\mathbb{N}^{*}} $ 
be a besselian Schauder frame of $E$. 
Let $\left( y_{n}^{*} \right) _{n\in\mathbb{N}^{*}} $
 be a sequence of elements of $E^{*} $ for which there exists a constant
  $\alpha\in [0,1[$ such that the following condition holds for each 
  $x\in E$ and $y^{*} \in E^{*}$:
 \begin{align*}
\underset{n= 1}{\overset{\infty}{\sum }} \left \vert (y^{*}_{n}-b^{*}_{n})(x)y^{*}(a_{n})\right \vert 
\leq \alpha \left\Vert x\right \Vert_{E}\left\Vert y^{*}\right \Vert_{E^{*}}
\end{align*} 
Then there exists a sequence $ \left(z_{n}^{*}\right)_{n\in\mathbb{N}^{*}} $  
 of elements of $E^{*}$ and a sequence $ \left(w_{n}\right)_{n\in\mathbb{N}^{*}} $  
  of elements of $E$ such that the following conditions holds:
\begin{enumerate}
\item $ \left(w_{n} \right)_{n\in\mathbb{N}^{*}}\equiv_{E}
 \left(a_{n} \right)_{n\in\mathbb{N}^{*}}  $. 
 \item $ \left(z_{n}^{*} \right)_{n\in\mathbb{N}^{*}}\equiv_{E^{*}}
 \left(y_{n}^{*} \right)_{n\in\mathbb{N}^{*}}$.
\item    $\left( \left(a_{n},z_{n}^{* }\right) \right) _{n\in\mathbb{N}^{*}} $
and $\left( \left(w_{n},b_{n}^{* }\right) \right) _{n\in\mathbb{N}^{*}} $
 are a besselian Schauder frames of $E$.
\end{enumerate}
\end{corollary}
$\textbf{Proof}$: Let us consider the paire 
$\left( \left(a_{n},y_{n}^{* }\right) \right) _{n\in\mathbb{N}^{*}} $  of $E$.
 We have then for each $x\in E$ and $y^{*}\in E^{*}$:
 \begin{align*}
\underset{n= 1}{\overset{\infty}{\sum }} \left \vert y^{*}_{n}(x)y^{*}(a_{n}-a_{n})+(y^{*}_{n}-b^{*}_{n})(x)y^{*}(a_{n})\right \vert 
\leq \alpha \left\Vert x\right \Vert_{E}\left\Vert y^{*}\right \Vert_{E^{*}}
\end{align*} 
It follows then, by virtue of the theorem \ref{thm.pert.b.Schauder.f}, that there exists a sequence
$ \left(z_{n}^{*}\right)_{n\in\mathbb{N}^{*}} $ of elements of $E^{*}$ and a sequence
$ \left(w_{n}\right)_{n\in\mathbb{N}^{*}} $ of elements of $E$ such that 
the conditions (1), (2) and (3) in the corollary are true. \qed
\begin{corollary} 
\label{perturbation.imp.dim}
Let $E$ be a weakly sequentially complete Banach space and 
$\left( \left(a_{n},b_{n}^{* }\right) \right) _{n\in\mathbb{N}^{*}} $ 
be a besselian Schauder frame of $E$. 
\begin{enumerate}
\item Assume that there exist $x^{0}_{1},...,x^{0}_{N}$ in $E$ $(N\in\mathbb{N}^{*})$
 and $\alpha\in [0,1[$ such that the following condition holds for each 
  $x\in E$ and $y^{*} \in E^{*}$:
  \begin{align*}
\underset{n= 1}{\overset{N}{\sum }} \left \vert b^{*}_{n}(x)y^{*}(x^{0}_{n}-a_{n})\right \vert 
+\underset{n= N+1}{\overset{\infty}{\sum }} \left \vert b^{*}_{n}(x)y^{*}(a_{n})\right \vert 
\leq \alpha \left\Vert x\right \Vert_{E}\left\Vert y^{*}\right \Vert_{E^{*}}
\end{align*}
Then $E$ is finite dimensional with $dim E \leq N$
\item Assume that there exist $y^{*}_{0,1},...,y^{*}_{0,N}$ in $E^{*}$ $(N\in\mathbb{N}^{*})$
 and $\alpha\in [0,1[$ such that the following condition holds for each 
  $x\in E$ and $y^{*} \in E^{*}$:
  \begin{align}
\label{condi.perturbation.imp.dim.2}
\underset{n= 1}{\overset{N}{\sum }} \left \vert (y^{*}_{0,n}-b^{*}_{n})(x)y^{*}(a_{n})\right \vert 
+\underset{n= N+1}{\overset{\infty}{\sum }} \left \vert b^{*}_{n}(x)y^{*}(a_{n})\right \vert 
\leq \alpha \left\Vert x\right \Vert_{E}\left\Vert y^{*}\right \Vert_{E^{*}}
\end{align}
Then $E$ is finite dimensional with $dim E \leq N$
\end{enumerate}
\end{corollary}
$\textbf{Proof}$: 
\begin{enumerate}
\item Let us consider the sequence
$  \left(x_{n} \right) _{n\in\mathbb{N}^{*}} $  of $E$. where:
\[x_{n}:=\left\{
\begin{array}{l}
x^{0}_{n},\;\;1\leq n \leq N \\
\;0,\;\; n\geq N+1
\end{array}
\right.\]   
 It follows from the corollary (\ref{pertubation.1}) that: 
 there exists a sequence $  \left(z_{n}^{*} \right) _{n\in\mathbb{N}^{*}} $ 
of elements of $E^{*}$ such that $\left( \left(x_{n},z_{n}^{* }\right) \right) _{n\in\mathbb{N}^{*}} $ is a besselian
 Schauder frame of $E$. Consequently, $E$ is finite dimentional with $dim E \leq N$.     
\item In view of the condition (\ref{condi.perturbation.imp.dim.2}) , let us consider the paire
$  \left( y_{n}^{* } \right) _{n\in\mathbb{N}^{*}} $ of $E^{*}$ where:
\[y^{*}_{n}:=\left\{
\begin{array}{l}
y^{*}_{0,n},\;\;1\leq n \leq N \\
\;0,\;\; n\geq N+1
\end{array}
\right.\]   
 It follows from the corollary (\ref{pertubation.2}) that: 
 there exists a sequence $  \left(w_{n}  \right) _{n\in\mathbb{N}^{*}} $ 
of elements of $E $ such that $\left( \left(w_{n},y_{n}^{* }\right) \right) _{n\in\mathbb{N}^{*}} $ is a besselian
 Schauder frame of $E$. Consequently, $E$ is finite dimentional with $dim E \leq N$. 
 \end{enumerate}  
\qed
\begin{remark}
We can obtain, by means of the corollary \ref{perturbation.imp.dim}, that if 
$\left( \left(a_{n},b_{n}^{* }\right) \right) _{n\in\mathbb{N}^{*}} $ of $E$  is a 
besselian Schauder frame of $E$ withe $\underset{n= 1}{\overset{+\infty}{\sum }} \left \Vert a_{n} \right \Vert_{E} \left \Vert b^{*}_{n} \right \Vert_{E^{*}}<+\infty$
 then $E$ is finite dimensional.
\end{remark}   
Indeed, assume that  $\underset{n= 1}{\overset{+\infty}{\sum }} \left \Vert a_{n} \right \Vert_{E} \left \Vert b^{*}_{n} \right \Vert_{E^{*}}<+\infty$,
 then there exists $N\in\mathbb{N}^{*}$ such that: 
 $\underset{n= N+1}{\overset{+\infty}{\sum }} \left \Vert a_{n} \right \Vert_{E} \left \Vert b^{*}_{n} \right \Vert_{E^{*}}<1$.
  Hence we have for each $x\in E$ and $y^{*}\in E^{*}$:
\begin{align*}
\underset{n= 1}{\overset{N}{\sum }} \left \vert b^{*}_{n}(x)y^{*}(a_{n}-a_{n})\right \vert 
+\underset{n= N+1}{\overset{\infty}{\sum }} \left \vert b^{*}_{n}(x)y^{*}(a_{n}-0)\right \vert 
\leq \alpha \left\Vert x\right \Vert_{E}\left\Vert y^{*}\right \Vert_{E^{*}}
\end{align*}
It follows that, by virtue of the corollary \ref{perturbation.imp.dim}, that $E$ is 
finite dimensional.\qed 
\begin{proposition}
Let $E$ be a Banach space (resp. a weakly sequentially complete Banach space) 
and $\left( \left(a_{n},b_{n}^{* }\right) \right) _{n\in\mathbb{N}^{*}} $ be a 
Schauder frame of $E$ such that $a_{n}\neq 0$ and $b^{*}_{n}\neq 0$ for each 
$n\in\mathbb{N}^{*}$ (resp. a besselian Schauder frame of $E$). Let
$V$ be a closed subspace of $E$ with $V \neq \lbrace 0\rbrace$ and $W$
 be a separable closed subspace of $E^{*}$ with 
 $W \neq \lbrace 0\rbrace$. Let $\mathcal{I}$ be an infinite subset of 
 $\mathbb{N}^{*}$ such that $\mathbb{N}^{*}\backslash \mathcal{I}$ is infinite.\\
 There exists a Schauder frame of $E$ (resp. a besselian Schauder frame of $E$)
 $ \left( \left(u_{n,1},z_{n}^{* }\right) \right) _{n\in\mathbb{N}^{*}}$
  and $ \left( \left(w_{n},v_{n,1}^{* }\right) \right) _{n\in\mathbb{N}^{*}}$
  such that:
\begin{enumerate}
\item $\underset{n\in \mathcal{I}}{\overline{span}}(u_{n,1})=V$ and 
$\underset{n\in \mathcal{I}}{\overline{span}}(z^{*}_{n})$ is isomorphic to $W$;
\item $\underset{n\in \mathcal{I}}{\overline{span}}(v^{*}_{n,1})=W$ and 
$\underset{n\in \mathcal{I}}{\overline{span}}(w_{n})$ 
is isomorphic to $V$;
\end{enumerate}
\end{proposition}
$\textbf{Proof}$: We set 
\[\mathcal{I}=\lbrace m_{1},m_{2},...\rbrace \;\; where,\;\;  m_{1}<m_{2}<...\]
\[\mathbb{N}^{*}\backslash \mathcal{I}=\lbrace n_{1},n_{2},...\rbrace \;\; where,\;\;  n_{1}<n_{2}<...\]
We consider the bijections $\sigma_{0}: \mathcal{I} \longrightarrow \mathbb{N}^{*} $
and $\sigma_{1}: \mathbb{N}^{*}\backslash \mathcal{I} \longrightarrow \mathbb{N}^{*} $
defined, for every $k\in \mathbb{N}^{*}$, by the formulas:
\[\sigma_{0}(m_{k})=k \;\;and\;\;\sigma_{1}(n_{k})=k \]
$E$ is separable, so $V$ is separable. Hence there exists a sequence 
$(c_{n})_{n\in \mathbb{N}^{*}}$ of elements of $E\backslash \lbrace 0\rbrace$ such
that: 
$$V=\underset{n\in \mathbb{N}^{*}}{\overline{span}}(c_{n})$$
Since $W$ is separable, it exists a sequence $(d_{n})_{n\in \mathbb{N}^{*}}$ of
elements of $E^{*}\backslash \lbrace 0\rbrace$ such that:
$$W=\underset{n\in \mathbb{N}^{*}}{\overline{span}}(d_{n}^{*})$$
Let us consider the paire $\mathcal{F}_{1}=\left( \left(u_{n},v_{n}^{* }\right) \right) _{n\in\mathbb{N}^{*}}$
 of $E$ defined by the formulas:
\[u_{n}:=\left\{
\begin{array}{l}
 a_{\sigma_{1}(n)},\;\;n\in \mathbb{N}^{*}\backslash  \mathcal{I}\\
 0, \;\;n\in \mathcal{I}
\end{array}
\right.
\;\;and \;\;\;\;v^{*}_{n}:=\left\{
\begin{array}{l}
b^{*}_{\sigma_{1}(n)},\;\;n\in \mathbb{N}^{*}\backslash  \mathcal{I}\\
 d^{*}_{\sigma_{0}(n)}, \;\;n\in \mathcal{I}
\end{array}
\right.\]
\begin{enumerate}
\item \label{Dem.subs.frame.1} We assume first that $\left( \left(a_{n},b_{n}^{* }\right) \right) _{n\in\mathbb{N}^{*}} $
 is a Schauder frame of $E$. It is then clear that the paires 
 $\mathcal{F}_{1}   $ is a Schauder frame of $E$. Let $(t_{n})_{n\in\mathbb{N}^{*}}$
 be a sequence of elements of $\mathbb{K}\backslash \lbrace 0 \rbrace$.
 We consider the paire  
 $\tilde{\mathcal{F}} =\left( \left(u_{n,1},v_{n,1}^{* }\right) \right) _{n\in\mathbb{N}^{*}}$
of $E$ defined by:
 \[u_{n,1}:=\left\{
\begin{array}{l}
 a_{\sigma_{1}(n)},\;\;n\in \mathbb{N}^{*}\backslash  \mathcal{I}\\
 t_{n}c_{\sigma_{0}(n)}, \;\;n\in \mathcal{I}
\end{array}
\right.
\;\;and \;\;\;\;
v^{*}_{n,1}:=\left\{
\begin{array}{l}
b^{*}_{\sigma_{1}(n)},\;\;n\in \mathbb{N}^{*}\backslash  \mathcal{I}\\
 d^{*}_{\sigma_{0}(n)}, \;\;n\in \mathcal{I}
\end{array}
\right.\]
We have:
\begin{align*}
 				 	 &\underset{n= 1}{\overset{\infty}{\sum }} \left \Vert v^{*}_{n,1}-v^{*}_{n}\right \Vert_{E^{*}}\left \Vert u_{n,1}\right \Vert_{E}
+\underset{u_{n}\neq 0}{\overset{}{\sum }}2K_{\mathcal{F}_{1}}\dfrac{ \left\Vert u_{n,1}-u_{n}\right \Vert_{E}}{\left \Vert u_{n}\right \Vert_{E}} 
+\underset{u_{n}= 0}{\overset{}{\sum }} \left \Vert v^{*}_{n}\right \Vert_{E^{*}}\left \Vert u_{n,1}\right \Vert_{E}
 \\
 				   	 &=\underset{n= 1}{\overset{\infty}{\sum }} \left \Vert v^{*}_{n}-v^{*}_{n}\right \Vert_{E^{*}}\left \Vert u_{n,1}\right \Vert_{E}
+\underset{n\in \mathbb{N}^{*}\backslash  \mathcal{I}}{\overset{}{\sum }}2K_{\mathcal{F}_{1}}\dfrac{ \left\Vert a_{\sigma_{1}(n)}-a_{\sigma_{1}(n)}\right \Vert_{E}}{\left \Vert a_{\sigma_{1}(n)}\right \Vert_{E}} 
+\underset{n\in\mathcal{I}}{\overset{}{\sum }} \left \Vert d^{*}_{\sigma_{0}(n)}\right \Vert_{E^{*}}\left \Vert t_{n}c_{\sigma_{0}(n)}\right \Vert_{E}\\
&=\underset{n\in\mathcal{I}}{\overset{}{\sum }} \left \vert t_{n} \right \vert  \left \Vert d^{*}_{\sigma_{0}(n)}\right \Vert_{E^{*}}\left \Vert c_{\sigma_{0}(n)}\right \Vert_{E}
				\end{align*}
Let us take $t_{n}\in \mathbb{K}\backslash \lbrace 0 \rbrace$, $ n\in \mathbb{N}^{*}$, in a way such that:
\[\underset{n\in\mathcal{I}}{\overset{}{\sum }} \left \vert t_{n} \right \vert  \left \Vert d^{*}_{\sigma_{0}(n)}\right \Vert_{E^{*}}\left \Vert c_{\sigma_{0}(n)}\right \Vert_{E}
<1\]
Hence, according to  theorem \ref{thm.pert.Schauder.f}, there exists a sequence $(z^{*}_{n})_{n\in \mathbb{N}^{*}}$ of elements
of $E^{*}$ and a sequence $(w_{n})_{n\in \mathbb{N}^{*}}$ of elements
of $E $ such that:
\begin{enumerate}
\item $ \left(w_{n} \right)_{n\in\mathbb{N}^{*}}\equiv_{E}
 \left(u_{n,1} \right)_{n\in\mathbb{N}^{*}}  $. 
 \item $ \left(z_{n}^{*} \right)_{n\in\mathbb{N}^{*}}\equiv_{E^{*}}
 \left(v_{n,1}^{*} \right)_{n\in\mathbb{N}^{*}}$.
\item    $\left( \left(u_{n,1},z_{n}^{* }\right) \right) _{n\in\mathbb{N}^{*}} $
and $\left( \left(w_{n},v_{n,1}^{* }\right) \right) _{n\in\mathbb{N}^{*}} $
 are a  Schauder frames of $E$.\\
  It follows from (a) and (b) that:
  \item 
  \begin{enumerate}
\item 
\begin{align*}
\underset{n\in \mathcal{I}}{\overline{span}}(u_{n,1}) 
&=\underset{n\in \mathcal{I}}{\overline{span}}\left (t_{n}c_{\sigma_{0}(n)}\right )  \\
 				   	 &=\underset{n\in \mathcal{I}}{\overline{span}}\left (c_{\sigma_{0}(n)}\right )   \\
  					 &=\underset{n\in \mathbb{N}^{*}}{\overline{span}}(c_{n}) \\
				&=V
 \end{align*}
 \item  $\left(z_{n}^{*} \right)_{n\in\mathcal{I}} \equiv_{E^{*}}
    \left(d_{\sigma_{0}(n)}^{*} \right)_{n\in\mathcal{I}}$. 
 So $\underset{n\in \mathcal{I}}{\overline{span}}\left (d_{\sigma_{0}(n)}^{*}\right )$ is isomorphic to 
 $\underset{n\in \mathcal{I}}{\overline{span}}(z^{*}_{n})$. But 
 $\underset{n\in \mathcal{I}}{\overline{span}}\left (d_{\sigma_{0}(n)}^{*}\right )
 =\underset{n\in \mathbb{N}^{*}}{\overline{span}}(d_{n}^{*})=W$. It follows that 
 $\underset{n\in \mathcal{I}}{\overline{span}}\left (z_{n}^{*}\right )$
 is isomorphic to $W$.
\end{enumerate}
\item
\begin{enumerate}
\item $ \left(w_{n} \right)_{n\in\mathcal{I}}\equiv_{E}
 \left(u_{n,1} \right)_{n\in\mathcal{I}}  $ . So $\underset{n\in \mathcal{I}}{\overline{span}}\left (w_{n}\right )$ is isomorphic to 
 $\underset{n\in \mathcal{I}}{\overline{span}}\left (u_{n,1}\right )=V$. On
 the other hand:  $\underset{n\in \mathcal{I}}{\overline{span}}\left (v^{*}_{n,1}\right )
 =\underset{n\in \mathcal{I}}{\overline{span}}\left (d_{\sigma_{0}(n)}^{*}\right )$
is isomorphic to $W$.
\end{enumerate}
\end{enumerate}
\item We assume first that $\left( \left(a_{n},b_{n}^{* }\right) \right) _{n\in\mathbb{N}^{*}} $
 is a besselian Schauder frame of $E$.
 It is then clear that the paire 
 $\mathcal{F}_{1}   $ is a besselian Schauder frame of $E$. We have: 
 \begin{align*}
 				 	 &\underset{n= 1}{\overset{\infty}{\sum }} 
\left \Vert v^{*}_{n,1}\right \Vert_{E^{*}}\left \Vert u_{n,1}-u_{n}\right \Vert_{E}
 +  \left \Vert v^{*}_{n,1}-v^{*}_{n}\right \Vert_{E^{*}}\left \Vert u_{n}\right \Vert_{E}\\
 				   	 &=\underset{n= 1}{\overset{\infty}{\sum }} 
\left \Vert v^{*}_{n,1}\right \Vert_{E^{*}}\left \Vert u_{n,1}-u_{n}\right \Vert_{E}
+  \left \Vert v^{*}_{n}-v^{*}_{n}\right \Vert_{E^{*}}\left \Vert u_{n}\right \Vert_{E}\\
&=\underset{n\in\mathcal{I}}{\overset{}{\sum }}   \left \Vert d^{*}_{\sigma_{0}(n)}\right \Vert_{E^{*}}\left \Vert t_{n}c_{\sigma_{0}(n)}\right \Vert_{E}\\
&=\underset{n\in\mathcal{I}}{\overset{}{\sum }} \left \vert t_{n} \right \vert  \left \Vert d^{*}_{\sigma_{0}(n)}\right \Vert_{E^{*}}\left \Vert c_{\sigma_{0}(n)}\right \Vert_{E}
				\end{align*}
Let us take $t_{n}\in \mathbb{K}\backslash \lbrace 0 \rbrace$, $ n\in \mathbb{N}^{*}$, in a way such that:
\[\underset{n\in\mathcal{I}}{\overset{}{\sum }} \left \vert t_{n} \right \vert  \left \Vert d^{*}_{\sigma_{0}(n)}\right \Vert_{E^{*}}\left \Vert c_{\sigma_{0}(n)}\right \Vert_{E}
<1\]
Hence, according to the corollary \ref{pert.bsf}, there exists a sequence $(z^{*}_{n})_{n\in \mathbb{N}^{*}}$ of elements
of $E^{*}$ and a sequence $(w_{n})_{n\in \mathbb{N}^{*}}$ of elements
of $E$ such that:
\begin{enumerate}
\item $ \left(w_{n} \right)_{n\in\mathbb{N}^{*}}\equiv_{E}
 \left(u_{n,1} \right)_{n\in\mathbb{N}^{*}}  $. 
 \item $ \left(z_{n}^{*} \right)_{n\in\mathbb{N}^{*}}\equiv_{E^{*}}
 \left(v_{n,1}^{*} \right)_{n\in\mathbb{N}^{*}}$.
\item    $\left( \left(u_{n,1},z_{n}^{* }\right) \right) _{n\in\mathbb{N}^{*}} $
and $\left( \left(w_{n},v_{n,1}^{* }\right) \right) _{n\in\mathbb{N}^{*}} $
 are a  besselian Schauder frames of $E$.\\
\end{enumerate} 
We prove in a similar way as in (\ref{Dem.subs.frame.1}) that:
$\underset{n\in \mathcal{I}}{\overline{span}}(u_{n,1})=V $
and $\underset{n\in \mathcal{I}}{\overline{span}}(z^{*}_{n}) $ is 
isomorphic to $W$; $\underset{n\in \mathcal{I}}{\overline{span}}(v^{*}_{n,1})=W$
and $\underset{n\in \mathcal{I}}{\overline{span}}(w_{n}) $ is 
isomorphic to $V$.\qed
\end{enumerate}


\end{document}